\documentclass[12pt,a4paper]{article}
\usepackage{hyperref,amsmath,amssymb}
\usepackage{graphicx}
\usepackage{geometry}
\usepackage{parskip}

\author{Rowena Ball\footnote{Professor Rowena Ball (email Rowena.Ball@anu.edu.au), Mathematical Sciences Institute,  Australian National University, Ngunnawal and Ngambri Lands, Canberra, ACT 2601 } 
} 

\date{}

\begin{document}

\title{Indigenous Mathematics \\
I. Smoke Telegraphy}

\setlength{\parskip}{6pt}
\maketitle

\begin{quotation}\noindent \textit{
The following article is the first in an occasional series for the Australian Mathematical Society Gazette on diverse aspects and topics of Indigenous mathematical knowledge. This is an important, but neglected, part of the mathematical heritage of humankind, and as such is the concern of the mathematics community as a whole. It is hoped that this and future articles may help to inspire mathematics researchers, students, and educators at tertiary and school levels who are seeking to widen  their mathematical horizons and develop course and research materials of broad cultural relevance.}
\end{quotation}
\begin{quotation}
\noindent \textit{
I would like to honour the Mithaka peoples of the Kurrawoolben and Kirrenderri (Diamantina) and Nooroondinna (Georgina) river channel country of SW Qld. The material in this article does not involve culturally restricted knowledge or images, and is shared with respect for the Mithaka ancestors and their descendants.}
\end{quotation}

\section*{}
\begin{figure}[h]
\centerline{\includegraphics[scale=1]{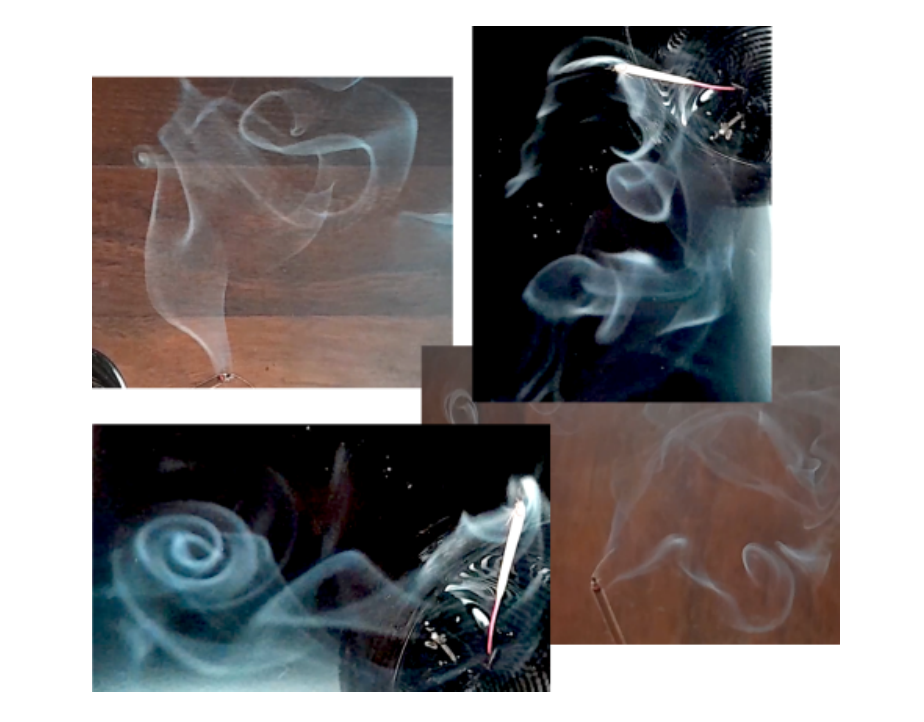}}
\caption{\label{figure1} The vortical structure of incense smoke flows is evident in this pastiche of stills from several experimental, but admittedly very amateur, home movies. No doubt you and your students could do much better!}
\end{figure}

Light a stick of incense. You will see a laminar plume of smoke emanating from the tip, which bifurcates into twin streams of counter-rotating vortices due to the shear in the flow giving rise to what is known as the Kelvin-Helmholtz instability, thence into other intriguing and beautiful secondary  structures, until the energy cascade dissipates them at the small scales (figure \ref{figure1}).  The vortices are formed as left-handed and right-handed pairs, a symmetry property (namely, chirality, which means they are non-superimposable mirror images) that may be used to encode and propagate information. 

And carry messages they did, in a long-distance smoke signalling technology that was practised by peoples of the Mithaka and other nations of the southwest Qld river channel country, and, indeed, throughout Australia \cite{Kerkhove:2021}. 

Around the turn of the 19th century, and for thousands of years previously, Mithaka lands and those of neighbouring nations were well-populated countries, a vast area rich in natural resources and supporting semi-permanent villages and lively trade and cultural exchanges, through which the well-traversed route network formed the greatest `silk road' of the continent \cite{Westaway:2021}. On this country most clans and individuals still lived their traditional ways of life, due in no small part to the terms of the Debney Peace Treaty of 1889 \cite{Griffith:2022}. Some people also worked on the cattle stations. 

Much of what we know about smoke telegraphy on Mithaka lands comes from the memoirs of Alice Duncan-Kemp \cite{Duncan-Kemp:1934,Duncan-Kemp:1952}, who grew up with the tribes of the district in the late 19th  and early 20th centuries on a cattle station called Mooraberrie. And here, in a village of native huts near the homestead, Alice introduces us to her childhood companion and carer, Mary Ann Coomindah, also known as Maghroolara, meaning singer, who was born on country around 1865.

Mary Ann Coomindah has an entry in the Australian Dictionary of Biography \cite{Coomindah:2022}, where her occupation is described as domestic servant and Indigenous healer or \textit{gdanja}. However, she was a true polymath who fulfilled many other roles involving advanced tribal knowledge.  Mary Ann belonged to the Emu tribe and was a leader of her people. She was highly educated in tribal knowledge, having passed her higher women’s Bora degrees. She was a teacher, and educated the young people of the tribes, girls and boys, on country. Mary Ann Coomindah was also an expert in smoke telegraphy or signalling, and I have no hesitation in calling her a practical  applied mathematician. 

Many of the details of smoke telegraphy have been lost or are sleeping, but what we \textit{can} say is that the technology took years of specialist, intergenerational, training to master, the code  was precise, necessarily mathematical, and cross-lingual, and the signals were interpreted, re-coded, and relayed on with perfect fidelity by recipients over vast distances. Alice Duncan-Kemp wrote at length about engineering the smoke signals in her memoirs. Although her writing style is idiosyncratic, the following excerpts leave no doubt about Mary Ann's applied mathematical knowledge, within cultural context, and the sophistication of the technology. The bold font emphases are mine. 
\begin{quotation}   \sffamily 
 
\noindent `Wreaths of smoke coiled themselves into \textbf{hoops} [rings] and twisted in \textbf{spirals} as messengers transmitted their talk in \textbf{smoke telegraphy} from tribe to tribe up and down the rivers. It was all part of a \textbf{vast network} of an \textbf{Aboriginal communication system}.’

\setlength{\parskip}{2mm}

\noindent ‘Other tribes \textbf{relayed} to still farther back tribes the \textbf{news} of the moment.' 

\noindent `Mary Ann and Bogie tapped out a message to the Walker-Di camp some 12 miles distant, now short, now long, and short again, \textbf{regulating} the waves until they became thin columns and smoke breaks [pulses], or delicate \textbf{spirals} with a \textbf{left to right curl} meaning ``white man'', instead of the more conventional \textbf{right to left spiral}.'

\noindent `Up flared the signal, \textbf{now light, now dark}, the spirals wavered in the breeze, then in short dense \textbf{pulses} alternating with dark columns\ldots an important message was being transmitted.'

\noindent `The code for smoke signalling is \textbf{intricate and very elaborate}. At a glance Mary Ann can tell whether a white man is mentioned in the faint smoke fluttering on the horizon, and who the man is and where he is.'

\noindent `A deep \textbf{black-winged column} of smoke shot up beside a \textbf{pale grey streak}. It was a \textit{Koto-unje}, a death signal'.

\noindent `A black smoke cloud twisted and twirled, now heavily in pulses, now lightly in pearly columns. For half an hour we watched the smoke ``talk'' from the \textit{Kibulyo} camp, sending out a  code message to be \textbf{relayed} by surrounding tribes.'
\end{quotation} 

Who would not want to know the news of the day? To send details of flood waters on the way, from rains in distant ranges, is an essential public service. When the dreaded black-winged column goes up, who will be invited to the funeral wake and ceremony (or perhaps more interestingly, who will not\ldots)?  

A widely-held belief, often deeply internalized, is  that mathematics is primarily of European (or Western) provenance, whence it is now firmly fixed in the curriculums of the whole world. Therefore, some may ask the question: Is there \textit{really} mathematical knowledge involved  in smoke telegraphy, which clearly is not part of the European or Western corpus? 

Well,  mathematics is nothing if not a body of knowledge that is culturally integrated into a  \textit{system}, and to exemplify this quality what could be more systematic than a communications technology? 
From the quotes above we note, or can infer,  the following mathematical aspects:

(a) Recognition and classification of three-dimensional symmetries (spirals of opposite handedness) and other shapes.

(b) Coding  that transcends language, using these symmetries and other structures such as rings, pulses with various time and length scales, enumeration of pulses, spin and colour. The relatively simple Morse code was used in electrical telegraphy --- we can only imagine the complexity of smoke codes, which must have been taught \textit{systematically}. A very limited collection of smoke codes is presented in \cite{Kerkhove:2021}.  

(c) Relaying, networking, and error correction.  

(d) The technical, spatial, and numerical  skills mastered in order to create and maintain smoke streams of a single chirality and colour or numerical sequences. 

Around that period in Western mathematics, classification of, and operations with, symmetries was a new field, and the operational vector and tensor calculus that enables fluid dynamicists to define and notate vorticity as the curl of the velocity field, in order  to conveniently model and compute vortical flows (e.g., by the Bureau of Meteorology), was only just being developed, largely by Heaviside (who, notably, was a telegraphist),  and in fact was highly controversial (especially in its notation). But Aboriginal women and men specialists, personified here by Mary Ann Coomindah, had already developed vortical smoke flows into a long-distance signalling technology for hundreds  of years or more.  

However, communication technologies typically are superseded. In 1872 the overland telegraph from Port Darwin to Port Augusta was completed, connecting Australia to the world via a submarine cable from Java. News of its progress had been smoke signalled by the tribes along the route! Electrical telegraphy was already being rolled out in the cities, by the end of that century telephony was well-established in urban areas, and from 1905 radio telegraphy was used too. Smoke signalling --- labour-intensive and restricted by weather conditions and to daylight hours --- became obsolete. But before that, it is fair to say that Aboriginal smoke telegraphy was \textit{the finest, most advanced, long-distance communication system in the world}. To create and send and decode and relay these signals, over vast distances, you \textit{have} to be a skilled practical mathematician, coder, and fluid dynamicist. And for that you have to have a systematic base of intercultural education and training passed down the generations and across cultures with different languages.

This article suggests various possible interesting and fun investigations for students. Can you, for example, develop and send a smoke code, at least in simulation? How, exactly, \textit{did} Mary Ann generate or separate right- and left-handed vortex streams? (Confession: my students and I made a complete mess of a single attempt! Persevere!) 

More importantly, though, research on smoke flows and composition has acquired new relevance and urgency, as global warming effects  create vastly enhanced patterns of wildland and farmland fires, and even urban fires. Indigenous-led research in this area is likely to yield transformational advances  in fire and landscape management on a warming Earth. 

There is more to Mary Ann Coomindah. I think she should be recognised as a true hero of Australian mathematics, of whom we all can be proud. She had other mathematical skills, namely, calendrical. Calendars are, or have been, kept by all cultures that  we know of, and are necessarily mathematical. Of all calendrical systems, the Gregorian calendar is possibly one of the worst. A discussion of \textit{that} claim, and of Indigenous calendars, will be elaborated in a future article.


The author acknowledges support from ARC grant IN230100053
\end{document}